\documentclass[12pt]{article}
\usepackage{amsxtra,amssymb,amsthm,amsmath,latexsym}

\theoremstyle{plain}
\newtheorem{theorem}{Theorem}[section]

\newtheorem{remark}[theorem]{Remark}
\newtheorem*{a1}{Assumption $A_1$}
\newtheorem*{a2}{Assumption $A_2$}
\newtheorem*{ad}{Assumption $A_d$}

\numberwithin{equation}{section}

\newcommand{\refT}[1]{Theorem~\ref{T:#1}}

\newcommand{\refR}[1]{Remark~\ref{R:#1}}

\def\R{{\mathbb R}}

\def\C{{\mathbb C}}

\def\dotf{\dot{f}}
\def\dotu{\dot{u}}
\def\dotA{\dot{A}}
\def\dotb{\dot{b}}
\def\dotB{\dot{B}}

\def\oH1{{\overset{\circ}{H}\kern-.04in{}^1}}

\def\bee{\begin{equation*}}
\def\eee{\end{equation*}}
\def\be{\begin{equation}}
\def\ee{\end{equation}}

\begin{document} 
\title{Continuity of
solutions to operator equations with respect to a parameter.  
}

\author{
A.G. Ramm\\
 Mathematics Department, Kansas State University, \\
 Manhattan, KS 66506-2602, USA\\
ramm@math.ksu.edu\\}

\date{}

\maketitle\thispagestyle{empty}

\begin{abstract}
\footnote{Math subject classification: 46T20, 47J05, 47J07 }
\footnote{key words: dependence on a parameter, operator equations  }

Let $A(k)u(k)=f(k) (1)$ be an operator equation, $X$ and $Y$ are
Banach spaces, $k\in\Delta\subset\C$ is a parameter, $A(k):X\to Y$ is
a map, possibly nonlinear. Sufficient conditions are given for
continuity of $u(k)$ with respect to $k$. Necessary and sufficient 
conditions are given for the continuity of $u(k)$ with respect to $k$
in the case of linear operators $A(k)$.  

\end{abstract}


\section{Introduction}\label{S:1}
Let $X$ and $Y$ be Banach spaces, $k\in\Delta\subset\C$ be a parameter, $\Delta$ be an open bounded set on a complex plane $\C$, 
$A(k):X\to Y$  be a map, possibly nonlinear, $f:=f(k)\in Y$ be a function.

Consider an equation \be\label{e1.1} A(k)u(k)=f(k). \ee We are interested
in conditions, sufficient for the continuity of $u(k)$ with respect to
$k\in\Delta$. There is a large literature (see eg. \cite{1}, \cite{2}) on
this subject. The novel points in our paper include necessary and
sufficient conditions for continuity of the solution to equation
\eqref{e1.1} and sufficient conditions for its continuity when the
operator $A(k)$ is nonlinear.

Consider separately the cases when $A(k)$ is a linear map and when $A(k)$
is a nonlinear map.

\begin{a1}\ $A(k):=X\to Y$ is a linear bounded operator, and 
 
a) equation \eqref{e1.1} 
is uniquely solvable for any $k\in\Delta_0:=\{k:|k-k_0|\leq r\}$, 
$k_0\in\Delta$, $\Delta_0\subset \Delta$, 

b)  $f(k)$ is continuous with respect to $k\in\Delta_0$,
$\sup_{k\in\Delta_0}\|f(k)\|\leq c_0$;

c) $\lim_{h\to 0} \sup_{\substack{k\in\Delta_0 \\ v\in M}} 
\|[A(k+h)-A(k)]v\|=0$, where $M\subset X$ is an arbitrary bounded set,

d) $\sup_{\substack{k\in\Delta_0 \\ f\in N}} \|A^{-1}(k)f\|\leq c_1$,
 where  $N\subset Y$ is an arbitrary bounded set, and $c_1$ may depend on 
$N$. \end{a1}

\begin{theorem} \label{T:1.1}
If Assumptions $A_1$ hold, then 
\be\label{e1.2}
\lim_{h\to 0} \|u(k+h)-u(k)\|=0.
\ee
\end{theorem}

\begin{proof}
One has 
\be\label{e1.3}
\begin{aligned}
u(k+h)-u(k)&=A^{-1}(k+h)f(k+h)-A^{-1}(k)f(k)\\
&=A^{-1}(k+h)f(k+h)-A^{-1}(k)f(k+h)+A^{-1}(k)f(k+h)-A^{-1}(k)f(k).
\end{aligned}\ee
\be\label{e1.4}
\|A^{-1}(k)[f(k+h)-f(k)]\|\leq c_1 \|f(k+h)-f(k)\|\to 0 \hbox{\quad 
as\quad} h\to 0.
\ee
\be\label{e1.5}
\begin{aligned}
\|A^{-1}(k+h)-A^{-1}(k) \|
   &= \|A^{-1}(k+h)[A(h+k)-A(k)]A^{-1}(k)\|\\
   &\qquad \leq c^2_1\|A(k+h)-A(k)\|\to 0 \hbox{\quad as\quad} h\to 0.
\end{aligned}
\ee 
From \eqref{e1.3}--\eqref{e1.5} and Assumptions $A_1$ the
conclusion of Theorem 1 follows. \end{proof}

\begin{remark}\label{R:1.2} Assumptions $A_1$ are not only sufficient for
the continuity of the solution to \eqref{e1.1}, but also necessary if one
requires the continuity of $u(k)$ uniform with respect to $f$ running 
through arbitrary bounded sets. Indeed, the necessity of the assumption a) 
is clear; that of the
assumption b) follows from the case $A(k)=I$, where $I$ is the identity 
operator; that of
the assumption c) follows from the case $A(k)=I$, $A(k+h)=2I$, 
$\forall h\not= 0$,
$f(k)=g\ \forall k\in\Delta_0$. Indeed, in this case assumption c) fails
and one has $u(k)=g$, $u(k+h)=\frac g 2$, so $||u(k+h)-u(k)||=\frac 
{||g||}2$ does not tend to zero as $h\to 0$.
 
To prove the necessity of the assumption d),
assume that $\sup_{k\in\Delta_0}\|A^{-1}(k)\|=\infty$. Then, by the
Banach-Steinhaus theorem, there is an element $f$ such that
$\sup_{k\in\Delta_0}\|A^{-1}(k)f\|=\infty$, so that
$\lim_{j\to\infty}\|A^{-1}(k_j)f\|=\infty$, $k_j\to k\in\Delta_0$. Then
$\|u(k_j)\| = \|A^{-1}(k_j)f\| \to\infty$, so $u(k_j)$ does not
converge to $u:=u(k)=A^{-1}(k)f$, although $k_j\to k$. \end{remark}

\begin{a2}\ $A(k):X\to Y$ is a nonlinear map and a), b), c) and d) of
Assumption $A_1$ hold, and the following assumption holds:

 e) $A^{-1}(k)$ {\it is a 
homeomorphism of $X$ onto $Y$ for each} $k\in\Delta_0$.
 \end{a2}

\begin{remark}\label{R1.3} Assumption e) is included in d) in the case of
a linear operator $A(k)$ because if $\|A(k)\|\leq c_2$ and assumption d)
holds, then 
$\|A^{-1}(k)\|\leq c_1$ and $A(k)$, $k\in \Delta_0,$ is an isomorphism of 
$X$ onto $Y$.
\end{remark}

\begin{theorem}\label{T:1.4}
If $A_2$ hold, then \eqref{e1.2} holds for the solution $u(k)$ to 
\eqref{e1.1}.
\end{theorem}

\begin{remark}\label{R:1.5}
Let us introduce the following assumption:

\begin{ad}: \ 

Assumptions $A_2$) hold and 

f) $\dotf(k):=\frac{df(k)}{dk}$ is continuous in $\Delta_0$,

g) $\dotA(u,k):= \frac{\partial A(u,k)}{\partial k}$ is continuous
with respect to (wrt)  $k$ in
$\Delta_0$ and wrt $u\in X$,

j) $\sup_{k\in\Delta_0} \|[A'(u,k)]^{-1}\|\leq c_3$,
where $A'(u,k)$ is the Fr\'echet derivative of  $A(u,k)$ and
$[A'(u,k)]^{-1}$ is continuous with respect to $u$ and $k$.
\end{ad}

{\bf Claim:} If Assumption $A_d$ holds, then 
\be\label{e1.6}
\lim_{h\to 0}\|\dotu(k+h)-\dotu(k)\|=0.
\ee
\end{remark}

\begin{remark}\label{R:1.6} If Assumptions $A_1$ hold except one: $A(k)$
is not necessarily a bounded linear operator, $A(k)$ may be unbounded,
closed, densely defined operator-function, then the conclusion of
\refT{1.1} still holds and its proof is the same. For example, let
$A(k)=L+B(k)$, where $B(k)$ is a bounded linear operator continuous with
respect to $k\in\Delta_0$, and $L$ is a closed, linear, densely defined
operator from $D(L)\subset X$ into $Y$. Then 
\bee
\|A(k+h)-A(k)\|=\|B(k+h)-B(k)\|\to 0 \hbox{\quad as\quad} h\to 0, \eee
although $A(k)$ and $A(k+h)$ are unbounded.  \end{remark}

In Section 2 proofs of \refT{1.4} and of \refR{1.5} are given.

\section{Proofs}

\begin{proof}[Proof of \refT{1.4}] One has: 
\bee
A(k+h)u(k+h)-A(k)u(k)=f(k+h)-f(k)=o(1) \hbox{\quad as\quad} h\to 0. 
\eee
Thus 
\bee A(k)u(k+h)-A(k)u(k)=o(1)-[A(k+h)u(k+h)-A(k)u(k+h)]. 
\eee 
Since
$\sup_{\{u(k+h): \|u(k+h)\|\leq c\}} \| 
A(k+h)u(k+h)-A(k)u(k+h)\|\underset{h\to
0}{\rightarrow}0$, one gets \be\label{e2.1} A(k)u(k+h)\to A(k)u(k)
\hbox{\quad as\quad} h\to 0. \ee 
By the Assumption $A_2$, item e), the operator $A(k)$ is a homeomorphism. 
Thus
\eqref{e2.1} implies \eqref{e1.2}.

 \refT{1.4} is proved.
\end{proof}

\begin{proof}[Proof of \refR{1.5}] First, assume that $A(k)$ is linear.
Then \be\label{e2.2} \frac{d}{dk} A^{-1}(k)=-A^{-1}(k)\dotA(k)A^{-1}(k),
\quad \dotA:=\frac{dA}{dk}. \ee Indeed, differentiate the identity
$A^{-1}(k)A(k)=I$ and get $\frac{dA^{-1}(k)}{dk}A(k)+A^{-1}(k)\dotA(k)=0$.
This implies \eqref{e2.2}. This argument proves also the existence of the
deriviative $\frac{dA^{-1}(k)}{dk}$. Formula $u(k)=A^{-1}(k)f(k)$ and the
continuity of $\dotf$ and of $\frac{dA^{-1}(k)}{dk}$ yield the existence
and continuity of $\dotu(k)$. \refR{1.5} is proved for linear operators
$A(k)$. \end{proof}

Assume now that $A(k)$ is nonlinear, $A(k)u:=A(u,k)$. Then one can
differentiate \eqref{e1.1} with respect to $k$ and get \be\label{e2.3}
\dotA(u,k)+A'(u,k)\dotu=\dotf, \ee where $A'$ is the Fr\'echet derivative
of $A(u,k)$ with respect to $u$. Formally one assumes that $\dotu$ exists,
when one writes \eqref{e2.3}, but in fact \eqref{e2.3} proves the
existence of $\dotu$, because $\dotf$ and $\dotA(u,k):=\frac{\partial
A(u,k)}{\partial k}$ exist by the Assumption $A_d$ and $[A'(u,k)]^{-1}$
exists and is an isomorphism by the Assumption $A_d$, item j). Thus,
\eqref{e2.3} implies 
\be\label{e2.4} \dotu=[A'(u,k)]^{-1}\dotf -
[A'(u,k)]^{-1}\dotA(u,k). \ee Formula \eqref{e2.4} and Assumption $A_d$
imply \eqref{e1.6}.

\refR{1.5} is proved. $\Box$

\section{Applications}\label{S:3}

\subsection{Fredholm equations depending on a parameter}\label{U:1}

Let 
\be\label{e3.1} Au:=u-\int_D b(x,y,k) u(y)dy:= [I-B(k)]u=f(k), \ee
where $D\subset R^n$ is a bounded domain, $b(x,y,k)$ is a function on
$D\times D\times\Delta_0$, $\Delta_0:=\{k: |k-k_0|<r\}$, $k_0>0$, 
$r>0$ is a
sufficiently small number. Assume that $A(k_0)$ is an isomorphism of
$H:=L^2(D)$ onto $H$, for example, $\int_D\int_D
|b(x,y,k_0)|^2dxdy<\infty$ and $N(I-B(k_0))=\{0\}$, where $N(A)$ is the
null-space of $A$. Then, $A(k_0)$ is an isomprohism of $H$ onto $H$ by the
Fredholm alternative, and Assumption $A_1$ hold if $f(k)$ is continuous
with respect to $k\in\Delta_0$ and \be\label{e3.2} \lim_{h\to 0}
\int_D\int_D |b(x,y,k+h)-b(x,y,k)|^2 dxdy=0 \qquad k\in\Delta_0.
\ee Condition \eqref{e3.2} implies that if $A(k_0)$ is an isomorphism of
$H$ onto $H$, then so is $A(k)$ for all $k\in\Delta_0$ if $|k-k_0|$ is
sufficiently small.

\refR{1.5} implies to \eqref{e3.1} if $\dotf$ is continuous with respect
to $k\in\Delta_0$, and $\dotb:=\frac{\partial b}{\partial k}$ is 
continuous
with respect to $k\in\Delta_0$ as an element of $L^2(D\times D)$.
Indeed, under these assumptions $\dotu=[I-B(k)]^{-1}(\dotf-\dotB(k)u)$ and
the right-hand side of this formula is continuous in $\Delta_0$.

\subsection{Semilinear elliptic problems}\label{U:2} 
Let \be\label{e3.3}
A_1(k)u:=Lu+g(u,k)=f_1(k), \ee where $L\geq m>0$ is an elliptic, second 
order,
self-adjoint, positive-definite operator with real-valued
coefficients in a bounded domain $D\subset\R^3$ with a smooth boundary,
and $g(u,k)$ is a smooth real-valued function on $\R\times\Delta_0$. Then
problem \eqref{e3.3} is equivalent to \eqref{e1.1} with
\be\label{e3.4} A(k)u:=u+L^{-1}g(u,k)=f(k):=L^{-1}f_1(k).
\ee
The operator $L^{-1}g(u,k)$ is compact in $C(D)$. Therefore
equation \eqref{e3.4} is solvable in $C(D)$ by the Schauder 
principle if 
the map $A(k)$ maps a ball $B(0,R):=B_R$ into itself for some $R>0$.
This happens if $g':=g'_u>0$ for $u>0$ and $\inf_{R>0} \frac {g(R)}R\leq 
m^{-1},$ 
where
$g(R):=\max_{k\in \Delta_0}|g(R,k)|$ and $||L^{-1}||_{C(D)\to C(D)}\leq 
m$.
Equation \eqref{e3.4}
has at most one solution if $g'>0$.
Assumptions $A_2$ can be verified, for example, if $g(u,k)$ is a 
smooth
function on $\R\times \Delta_0$ and 
$g'\geq 0$. In
this case $\|A^{-1}(k,f)\|\leq c\|f\|$, and $A^{-1}(k)$ is a continuous
operator defined on all of $H:=H^2_0(D)$, where $H$ is a real Hilbert 
space, 
for any fixed
$k\in\Delta_0$.
If, for example, $L=-\Delta +k^2$ is the Dirichlet operator in $D\subset 
\R^3$, then 
$L^{-1}$ is a positive-definite integral operator with the kernel
$0\leq G(x,y)<\frac {\exp(-k|x-y|)}{4\pi |x-y|}$, and $m\leq \frac 
{\int_0^{ka}e^{-s}sds}{k^2}$, where $a$ is the radius of $D$,
that is, $2a:=\sup_{x,y\in D} |x-y|$.


\begin{thebibliography}{1000} 


\bibitem[1]{1}
Kantorovich, L., Akilov, G., 
Functional analysis, Pergamon Press, New York, 1982

\bibitem[2]{2}
Ramm, A.G., 
Some theorems on equations with parameters
in Banach space, Doklady Acad. of Sci. Azerb. SSR, 22, (1966), 3-6. 

\end{thebibliography}
\end{document}